\documentclass[12pt]{article}
\usepackage[utf8]{inputenc}

\usepackage[textwidth=6.7in,top=1.2in,bottom=1.2in]{geometry}
\usepackage{amsthm,amsmath,amssymb,bbm,bm}
\usepackage{natbib}
\usepackage{multirow}
\usepackage[pdftex]{graphicx}
\usepackage{wrapfig}
\usepackage{tikz}
\usepackage{centernot}
\usepackage{url}
\usepackage{algorithmic}
\usepackage[linesnumbered, ruled, vlined]{algorithm2e}
\usepackage{dsfont}
\usepackage{titling}
\usepackage{relsize}
\usepackage{rotating}
\usepackage{enumitem}
\usepackage{titling}
\usepackage{float}
\usepackage{hyperref}
\usepackage{setspace}
\newtheorem{assumption}{Assumption}

\newtheorem{theorem}{Theorem}[section]

\newcommand{\sign}{{\text{sign}}}

\newcommand{\cC}{{\mathcal{C}}}

\newcommand{\cD}{{\mathcal{D}}}

\newcommand{\indep}{\perp \!\!\! \perp}

\title{On a necessary and sufficient identification condition of optimal treatment regimes with an instrumental variable}
\author{Yifan Cui and Eric Tchetgen Tchetgen\\ National University of Singapore, University of Pennsylvania}
\date{}

\begin{document}
\maketitle

\abstract{
Unmeasured confounding is a threat to causal inference and individualized decision making. Similar to \cite{cui2020,qiu2020,Han2020},
we consider the problem of identification of optimal individualized treatment regimes with a valid instrumental variable. 
 \cite{Han2020} provided an alternative identifying condition of optimal treatment regimes using the conditional Wald estimand of \cite{cui2020,qiu2020} when treatment assignment is subject to endogeneity and a valid binary instrumental variable  
is available.
In this note, we provide a necessary and sufficient condition for identification of optimal treatment regimes using the conditional Wald estimand.  Our novel condition is necessarily implied by those of  \cite{cui2020,qiu2020,Han2020} and may continue to hold in a variety of potential settings not covered by prior results.
}

\vspace{0.6cm}

\noindent {\bf keywords:}
 Individualized decision making, policy making, optimal treatment regimes, unmeasured confounding, sign identification, conditional average treatment effect

\section{Introduction}
Estimating optimal treatment regimes is a central task for precision medicine. In the health sciences and medicine, 
an individualized treatment regime provides a personalized treatment strategy for each
patient in the population based on individual characteristics.
A prevailing strand of work has been devoted to estimating optimal treatment regimes (\cite{robins2000marginalSM,murphy2003optimal,qian2011performance,zhao2012estimating,zhang2012estimating} and many others), we refer to \cite{chakraborty2013statistical,kosorok2019review,tsiatis2019dynamic} for an up-to-date literature review on dynamic treatment regimes.

Recently, there has been a fast-growing literature on estimating individualized treatment regimes based on observational
studies subject to potential unmeasured confounding \citep{kallus2018confounding,yadlowsky2018bounds,pmlr-v89-kallus19a,han2019optimal,cui2020,qiu2020,Han2020}.
In particular, \cite{cui2020} and \cite{qiu2020} tackled the problem of individualized decision making/estimating individualized treatment regimes by leveraging a valid instrumental variable (IV) to account for potential unmeasured confounding. 
\cite{Han2020} provided an alternative identifying condition under which the conditional Wald estimand on which identification is based on in \cite{cui2020,qiu2020} continues to identify optimal treatment regimes.

In this paper, we further relax sufficient identifying conditions considered by \cite{cui2020,qiu2020,Han2020}, and introduce the concept of identifying optimal treatment regimes by only identifying the sign of conditional average treatment effect (CATE). We propose a necessary and sufficient identification condition based on the possibility of identifying the sign of the CATE from the conditional Wald estimand without necessarily being able to identify the CATE nor the population value function for any given treatment regime. We illustrate the result by exploring realistic scenarios in which conditional optimal treatment regime can be identified in settings in which identifying conditions of  \cite{cui2020,qiu2020,Han2020} are not met. 
It is notable that our Theorem~\ref{thm:1} allows for identification of optimal treatment regimes even when an unmeasured confounding factor is an effect modifier of both the first stage association between the IV and the endogenous treatment, and of the causal effect of the endogenous treatment on the outcome; a possibility first recognized by \cite{Han2020} which we have hereby significantly expanded upon.

To conclude this section, we briefly introduce notation used throughout the paper. Let $Y$ denote the outcome of interest and $A \in \{-1,1\}$ be a binary treatment indicator.
Throughout it is assumed without loss of generality that larger values of $Y$ are more desirable.
Suppose that $U$ is an unmeasured confounder of the effect of $A$ on $Y$.  Suppose also that one has observed a pre-treatment binary IV $Z \in \{-1,1\}$. Let $L$ denote a set of fully observed pre-IV covariates.
Throughout we assume the complete data are independent and identically distributed realizations of $(Y, L, A, Z, U)$; thus the observed data are $(Y,L,A,Z)$.

\section{IV approaches to optimal treatment regimes}

We wish to identify an optimal treatment regime $\cD^*$, which is a mapping from the patient-level covariate space to the treatment
space $\{-1, 1\}$ that maximizes the corresponding expected potential outcome for the entire population, i.e.,
\begin{align*}
\cD^* =\arg\max_{\cD} E_L \left[  E_{Y_{\cD}} [Y_{\cD(L)}|L] \right],
\end{align*}
where $Y_a$  is a person's potential outcome under an intervention that sets the treatment to value $a$, and $Y_{\cD(L)}$ is the potential outcome under a hypothetical intervention that assigns treatment according to the regime $\cD$, i.e.,
\begin{align*}
Y_{\cD(L)} \equiv Y_{1}I\{\cD(L)=1\}+Y_{-1}I\{\cD(L)=-1\},
\end{align*}
where $I\{\cdot\}$ is the indicator function.
Optimal individualized treatment regimes can alternatively be written as
\begin{align}
\cD^*(L) = \sign\{ \gamma(L) \}, \label{eq:opt2}
\end{align}
where $\gamma(L) \equiv  E(Y_1-Y_{-1}|L)\neq 0$.
Throughout the paper, we make the standard consistency and positivity assumptions as in \cite{cui2020}.

A significant amount of work has been devoted to estimating optimal treatment regimes relying on the following unconfoundedness assumption:
\begin{assumption}\emph{(Unconfoundedness)}
$Y_a \indep A| L$ for $a=\pm 1$.
\label{asm:unconfoundedness}
\end{assumption}
The assumption essentially rules out the existence of an unmeasured factor $U$ that confounds the effect of $A$ on $Y$ upon conditioning on $L$. It is straightforward to verify that
under Assumption~\ref{asm:unconfoundedness}, one can identify the value function \citep{qian2011performance} $E[Y_{\cD(L)}]$ for a given treatment regime $\cD$.
Furthermore, optimal treatment regimes in Equation~\eqref{eq:opt2} are identified from the observed data
\begin{align*}
\cD^*(L) = \sign\{ E(Y|A=1,L) - E(Y|A=-1,L) \}.
\end{align*}
As established by \cite{qian2011performance}, learning optimal treatment regimes under Assumption~\ref{asm:unconfoundedness} can be formulated as
\begin{align*}
\cD^*=\arg\max_{\cD} E\left[\frac{YI\{\cD(L)=A\}}{f(A|L)}\right].
\end{align*}
\cite{zhang2012robust} proposed to directly maximize the value function over a parametrized set of functions.
Rather than maximizing the above value function, \cite{zhao2012estimating,zhang2012estimating,rubin2012statistical} transformed the above problem into a weighted classification approach, which was shown to have appealing robustness properties, particularly in a randomized study where no model assumption on $Y$ is needed.

Instead of relying on Assumption~\ref{asm:unconfoundedness}, we allow for unmeasured confounding.
Let $Y_{z,a}$ denote the potential outcome had, possibly contrary to fact, a person's IV and treatment value been set to $z$ and $a$, respectively.
 Suppose that the following assumption holds.
\begin{assumption}\emph{(Latent unconfoundedness)}
$Y_{z,a} \indep (Z, A)|L, U$ for $z,a = \pm 1$.
\label{asm:unconfoundedness2}
\end{assumption}
This assumption essentially states that together $U$ and $L$ would in principle suffice to account for any confounding bias.
Because $U$ is not observed, we propose to account for it when a valid IV $Z$ is available that satisfies the following standard IV assumptions \citep{angrist1996}:

\begin{assumption}\emph{(IV relevance)} $Z \centernot{\indep} A|L$.
\label{IV Relevance}
\end{assumption}

\begin{assumption}\emph{(Exclusion restriction)} $Y_{z,a}=Y_a$ for $z,a=\pm 1$ almost surely.
\label{Exclusion Restriction}
\end{assumption}

\begin{assumption}\emph{(IV independence)} $Z \indep  U |L$.
\label{IV Independence}
\end{assumption}

\begin{assumption}\emph{(IV positivity)} $0<f\left(  Z=1|L\right)<1$ almost surely.
\label{asm:IV positivity}
\end{assumption}

Assumptions~\ref{IV Relevance}-\ref{IV Independence} are well-known IV conditions, while Assumption~\ref{asm:IV positivity} is needed for nonparametric identification \citep{10.1093/ije/29.4.722,hernan2006epi}.
Assumption~\ref{IV Relevance} requires that the IV is associated with the treatment conditional on $L$. Note that Assumption~\ref{IV Relevance} does not rule out confounding of the $Z$-$A$ association by an unmeasured factor, however, if present, such factor must be independent of $U$.
 Assumption~\ref{Exclusion Restriction} states that there can be no direct causal effect of $Z$ on $Y$ not mediated by $A$.
 Assumption~\ref{IV Independence} states that the direct causal effect of $Z$ on $Y$ would be identified conditional on $L$ if one could intervene on $A=a$.
Hereinafter, we refer to Assumptions~\ref{IV Relevance}-\ref{asm:IV positivity} as core IV assumptions. Figure~\ref{fig:1} provides a graphical representation of
Assumptions \ref{Exclusion Restriction} and \ref{IV Independence}.

\tikzstyle{var} = [circle,  very thick,draw=black,fill=white,minimum size=10mm]
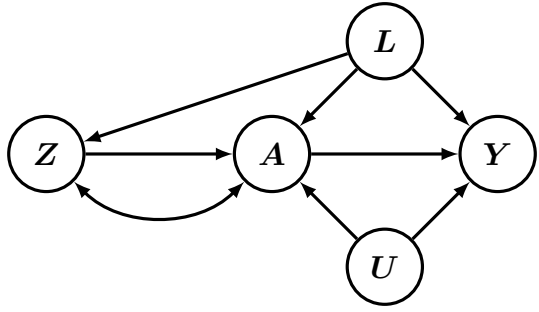
\begin{figure}[h]
\centering
\begin{tikzpicture}[>=latex]
    \node (A) [var, xshift=0cm,yshift=0cm]{$\boldsymbol{A}$};
    \node (U) [var, xshift=1.5cm,yshift=-1.5cm]{$\boldsymbol{U}$};
    \node (Y) [var, xshift=3cm,yshift=0cm]{$\boldsymbol{Y}$};
    \node (X) [var, xshift=1.5cm,yshift=1.5cm]{$\boldsymbol{L}$};
    \node (Z) [var, xshift=-3cm,yshift=0cm]{$\boldsymbol{Z}$}
    edge[<->, bend right=45, very thick](A)
    ;
    \draw[->, very thick] (X) -- (Z);
    \draw[->, very thick] (X) -- (A);
    \draw[->, very thick] (X) -- (Y);
    \draw[->, very thick] (Z) -- (A);
    \draw[->, very thick] (A) -- (Y);
    \draw[->, very thick] (U) -- (A);
    \draw[->, very thick] (U) -- (Y);
\end{tikzpicture}
\caption{A causal graph with unmeasured confounding. The bi-directed arrow between $Z$ and $A$ indicates the possibility that there may be unmeasured common causes confounding their association. \label{fig:1}}
\end{figure}

 These four core IV assumptions together do not suffice for point identification of the counterfactual mean and average treatment effect. \cite{cui2020} showed that 
 $\cD^*$ is nonparametrically identified by 
\begin{align}
\arg\max_{\cD}E \left[\frac{ZAYI\{A=\cD(L)\}}{\delta(L)f(Z|L) } \right], \label{f1}
\end{align}
and
\begin{align}
\arg \max_\cD  E \left[\frac{YI\{Z=\cD(L)\}}{\delta(L)f(Z|L) } \right], \label{f2}
\end{align}
under 
no unmeasured common effect modifier assumption, i.e.,
$$ Cov\left\{\widetilde \delta(L,U), \widetilde \gamma(L,U)|L\right\} = 0,$$ almost surely,
where $\widetilde \delta(L,U) \equiv \Pr(A=1|Z=1,L,U)-\Pr(A=1|Z=-1,L,U)$ and $\widetilde \gamma(L,U) \equiv E(Y_1-Y_{-1}|L,U)$, respectively,
or 
independent compliance type assumption, i.e., 
$$ \delta(L) \equiv \Pr(A=1 | Z=1,L)-\Pr(A=1| Z=-1,L)=\widetilde \delta(L,U)~ \text{almost surely}.$$
\cite{Han2020} proposed the following alternative identifying assumption for \eqref{f1} and \eqref{f2} to identify optimal treatment regimes given for a causal IV: The following two conditions hold given $L$,\\
(a) either $\widetilde \gamma(L,U) \geq 0$ or $\widetilde \gamma(L,U) \leq 0$ almost surely; and \\
(b) either $\widetilde \delta(L,U) \geq 0$ or $\widetilde \delta(L,U) \leq 0$ almost surely.

\section{Individualized treatment regimes: Identifying the sign of CATE}

We note that in order to identify optimal treatment regimes, one only needs to identify the sign of CATE.  
Therefore, in the following theorem, we provide sufficient and necessary conditions for Equations~\eqref{f1} and \eqref{f2} identifying $\cD^*$, or equivalently, the conditional Wald estimand
${\cC(L)}/{\delta(L)}$
in \cite{wald1940fitting,wang2018bounded} having the same sign as $\gamma(L)$, where $\cC(L)\equiv E(Y|Z=1,L) - E(Y|Z=-1,L)$.
\begin{theorem}\label{thm:1}
Under Assumptions~\ref{asm:unconfoundedness2}-\ref{asm:IV positivity}, 
the following condition is necessary and sufficient for  Equation~\eqref{f1} or \eqref{f2} identifying $\cD^*(L)$, 
\begin{align}\label{eq:eq1}
         E \bigg( \frac{\widetilde \gamma(L,U)}{\gamma(L)} \times \frac{\widetilde \delta(L,U)}{ \delta(L)}\bigg|L \bigg) > 0,
\end{align} 
or equivalently, 
\begin{align}\label{eq:eq2}
         Cov \bigg( \frac{\widetilde \gamma(L,U)}{\gamma(L)},\frac{\widetilde \delta(L,U)}{ \delta(L)}\bigg|L \bigg) > -1.
\end{align}
\end{theorem}

The necessary and sufficient identifying condition clearly holds if 
\begin{align}\label{eq:rt}
Cov\bigg( \frac{\widetilde \gamma(L,U)}{\gamma(L)},\frac{\widetilde \delta(L,U)}{ \delta(L)} \bigg|L\bigg) \geq 0,
\end{align}
which essentially states that $\frac{\widetilde \gamma(L,U)}{\gamma(L)}$ and $\frac{\widetilde \delta(L,U)}{ \delta(L)}$ are positively correlated conditional on $L$.
This would hold for instance if both functions are either non-increasing or non-decreasing as functions of $U$ whenever $U$ is scalar.  Importantly, this assumption may hold if either ${\widetilde \gamma(L,U)}$ or $\widetilde \delta(L,U)$ is positive for some values of $U$ and negative for other values conditional on $L$ such that the sign of say $\widetilde \delta(L,U)$ does not always agree with that of $\delta(L)$ therefore invalidating the assumption of \cite{Han2020} (even if a causal IV is assumed in \cite{Han2020}). 

The theorem establishes that optimal treatment regimes are identified by the sign of the conditional Wald estimand \citep{wald1940fitting,wang2018bounded} in the setting where individuals' decision to uptake the intervention is concordant with an anticipated benefit from the intervention. Specifically, consider the following definition of concordant treatment uptake with anticipated benefit from the intervention in the special case where $U$ can be taken as univariate:  
a) $E(Y_1-Y_{-1}|L,U)$ and $E(A|Z=1,L,U)-E(A|Z=-1,L,U)$ are both non-decreasing or non-increasing in $U$ so that as anticipated expected treatment benefit increases with $U$, expected treatment uptake likewise increases in $U$; b) $\sign\{ E(Y_1-Y_{-1}|L)\}=\sign \{E(A|Z=1,L)-E(A|Z=-1,L)\}$ so that conditional on $L$, patients' expected decision to uptake the intervention matches  their expected treatment benefit. 
We refer to subjects with an ability to fulfill both a) and b) as rational agents with perfect anticipation. It is then clear according to the theorem that if all subjects are rational agents with perfect anticipation, optimal treatment regimes can be identified on the basis of the sign of the conditional Wald estimand. Importantly, the theorem allows for some portion of the population of subjects failing to be rational possibly due to imperfect anticipation of expected treatment benefit; provided that they do not offset the contribution to the conditional Wald estimand from rational agents with perfect anticipation of expected treatment benefit.

Our necessary and sufficient condition is also clearly implied by those of \cite{cui2020,qiu2020,Han2020}. 
The assumption of  \cite{cui2020} corresponds to the case
$Cov\left\{\widetilde \delta(L,U), \widetilde \gamma(L,U)|L\right\} = 0;$
while that of  \cite{Han2020}  given for a causal IV: Conditional on $L$,\\
(a) either $\widetilde \gamma(L,U) \geq 0$ or $\widetilde \gamma(L,U) \leq 0$ almost surely; and \\
(b) either $\widetilde \delta(L,U) \geq 0$ or $\widetilde \delta(L,U) \leq 0$ almost surely;\\ 
imply Equation~\eqref{eq:eq1}
as it implies the sign of ${\widetilde \gamma(L,U)}$ agrees with the sign of ${\gamma(L)}$ almost surely and likewise the sign of $\widetilde \delta(L,U)$ agrees with that of  $\delta(L)$ almost surely.

\section{Three levels of assumptions for identifying optimal treatment regimes}

In this section, as a summary, one may categorize various assumptions for identifying optimal treatment regimes into three levels, namely depending on whether one can identify the value function, the CATE, or only the sign of CATE in Table~\ref{tb:1}.

In particular, as for identifying the CATE, Assumption~A5b(1b) of \cite{qiu2020} is a special case of Assumption~7 of \cite{cui2020},  while Assumption~A5b(1a) of \cite{qiu2020} relaxes IV independence to uncorrelated IV. 
As for identifying the sign of CATE, while both imply Equation~\eqref{eq:eq1} or \eqref{eq:eq2}, Assumption~A of \cite{Han2020} and Equation~\eqref{eq:rt} do not necessarily imply each other. 
Moreover, nonparametric IV bounds such as Balke-Pearl bounds \citep{Balke1997}
can be used to identify the sign of CATE when not covering zero \citep{cui2020}.

Interestingly, we point out that while, on the left side of Table~\ref{tb:1}, the quantity that is identified is stronger from bottom to top, identifying assumptions in the top do not necessarily imply those in the bottom on the right side.  
For instance, Assumption~7 of \cite{cui2020} and Assumption~A of \cite{Han2020} do not necessarily imply each other because no heterogeneity in $U$ of the compliance type or the average additive treatment effect on the outcome implies Assumption 7 of \cite{cui2020} without implying Assumption A of \cite{Han2020}; and the latter assumption does not necessarily imply the no common unmeasured effect modifier condition.

\begin{table}[H]
\centering
\resizebox{\columnwidth}{!}{
\begin{tabular}{|c|c|}
\hline
Quantity Identified & Identifying Assumptions    \\
\hline
 the value function & Assumption~8 of \cite{cui2020},  Assumption~A5b(2) of \cite{qiu2020}   \\\hline
 the CATE  & Assumption~7 of \cite{cui2020}, Assumption~A5b(1) of \cite{qiu2020}  \\ \hline
the sign of CATE &   Assumption~A of \cite{Han2020}, Equation~\eqref{eq:eq1} or \eqref{eq:eq2} or \eqref{eq:rt}, IV CATE bounds when not covering 0\\
\hline
\end{tabular}
}
\caption{A categorization of various identifying assumptions in \cite{cui2020}, \cite{qiu2020}, and \cite{Han2020} for identifying the value function, CATE, and sign of CATE. \label{tb:1}}
\end{table}

\newpage
\appendix

\begin{center}
{\LARGE \textbf{Appendix}}
\end{center}

\section{Proof of Theorem~\ref{thm:1}}
\begin{proof}
From the proof of Theorems~2.1 and 2.2 in \cite{cui2020}, we have  
\begin{align*}
&E \left[\frac{ZAYI\{A=\cD(L)\}}{\delta(L)f(Z|L) } \right]\\
=& E \left[ \frac{  \left[\Pr(A=1| Z=1,L,U)-\Pr(A=1| Z=-1,L,U) \right]I\{\cD(L)=1\}E[Y_1|L,U]  }{\delta(L)}  \right] \\
+ & E \left[ \frac{  \left[\Pr(A=1| Z=1,L,U)-\Pr(A=1| Z=-1,L,U) \right]I\{\cD(L)=-1\}E[Y_{-1}|L,U]  }{\delta(L)}  \right],
\end{align*}
and
\begin{align*}
 & E \left[\frac{YI\{Z=\cD(L)\}}{\delta(L)f(Z|L) } \right]\\
 = & E \left[ \frac{  \left[\Pr(A=1| Z=1,L,U)-\Pr(A=1| Z=-1,L,U) \right]I\{\cD(L)=1\}E[Y_1|L,U]  }{\delta(L)}  \right] \\
+ & E \left[ \frac{  \left[\Pr(A=1| Z=1,L,U)-\Pr(A=1| Z=-1,L,U) \right]I\{\cD(L)=-1\}E[Y_{-1}|L,U]  }{\delta(L)}  \right] \\
+ & E \left[ \frac{ \Pr(A=1| Z=-1,L,U) E[Y_1|L,U]}{ \delta(L)}  \right] \\
+ & E \left[ \frac{ \Pr(A=-1| Z=-1,L,U) E[Y_{-1}|L,U]}{ \delta(L)}  \right].
\end{align*}

It is easy to see that 
\begin{align*}
\arg\max_{\cD}E \left[\frac{ZAYI\{A=\cD(L)\}}{\delta(L)f(Z|L) } \right],
\end{align*}
and
\begin{align*}
\arg \max_\cD  E \left[\frac{YI\{Z=\cD(L)\}}{\delta(L)f(Z|L) } \right], 
\end{align*}
equal to 
\begin{align}\label{eq:f3}
\arg \max_\cD  E \left[\frac{\widetilde \gamma(L,U)\widetilde \delta(L,U)}{\delta(L)}I\{\cD(L)=1\} \right],
\end{align}
which is denoted by $\widetilde \cD$.
This completes the proof as $\widetilde \cD(L)$
necessarily agrees with the sign of $\gamma(L)$ whenever  Equation~\eqref{eq:eq1} holds.

The proof can also be conducted from the perspective of conditional Wald estimand.
Recall that 
\begin{align*}
\cC(L)
= & E \big \{E[Y | Z = 1,L,U]  -  E[Y | Z = -1, L,U] \big|L \big\} \\
= & E \big \{E[Y_1 \frac{1+A}{2}| Z = 1,L,U]  +  E[Y_{-1} \frac{1-A}{2}| Z = 1, L,U] \big|L \big\}  \\
& -  E \big \{E[Y_1 \frac{1+A}{2} | Z = -1,L,U]  +  E[Y_{-1}\frac{1-A}{2} | Z = -1, L,U]\big|L \big\} \\
= & E \bigg(E[Y_1 - Y_{-1} | L,U] \Big\{\Pr[A=1| Z = 1, L, U] - \Pr[A=1| Z = -1,  L,U]\Big\} \bigg|L \bigg)\\
= & E \Big( \widetilde \gamma(L,U) \widetilde \delta(L,U)  \Big|L \Big).
\end{align*}

Thus, we have that 
\begin{align*}
        \frac{\cC(L)}{\delta(L)}= E \Big( \frac{\widetilde \gamma(L,U) \widetilde \delta(L,U)}{\delta(L)}  \Big|L \Big).
\end{align*}
Subsequently, 
\begin{align*}
        \frac{\cC(L)}{\delta(L)}/\gamma(L) = E \bigg( \frac{\widetilde \gamma(L,U)}{\gamma(L)}\times \frac{\widetilde \delta(L,U)}{ \delta(L)}\bigg|L \bigg).
\end{align*}
Also note that 
\begin{align*}
     & Cov \bigg( \frac{\widetilde \gamma(L,U)}{\gamma(L)},\frac{\widetilde \delta(L,U)}{ \delta(L)}\bigg|L \bigg) +1 \\
     = & E \bigg( \frac{\widetilde \gamma(L,U)}{\gamma(L)} \times \frac{\widetilde \delta(L,U)}{ \delta(L)}\bigg|L \bigg),
\end{align*}
which completes the proof as the sign of $\gamma(L)$ will necessarily agree with that of the conditional Wald estimand whenever the expression in above display is positive.

\end{proof}

\section{Identifying assumptions in Table~\ref{tb:1}}

For easy referencing, we provide identifying assumptions appeared in Table~\ref{tb:1}.

\begin{itemize}
 \item Assumption~7 of \cite{cui2020} \emph{(No unmeasured common effect modifier)}:
 $Cov\left\{\widetilde \delta(L,U), \widetilde \gamma(L,U)|L\right\} = 0$ almost surely.

\item Assumption~8 of \cite{cui2020} \emph{(Independent compliance type)}:
$ \delta(L)=\widetilde \delta(L,U)$ almost surely.
\label{asm:strong}
 
\item Assumption~A5b(1) of \cite{qiu2020}:  Both conditions below hold:\\
(a) (Uncorrelated IV) $Cov(Y_{-1},Z |L)=0$ almost surely;\\
(b) (No unmeasured treatment-outcome effect modification)
$E[Y_1-Y_{-1}|L,U]=E[Y_1-Y_{-1}|L]$ almost surely.

\item Assumption~A5b(2) of \cite{qiu2020}:  Both conditions below hold:\\
(a) (Independent IV) $Z$ and $U$ are independent given $L$;\\
(b) (Independent compliance)
$E[A_Z|Z=1,L,U]-E[A_Z|Z=-1,L,U]=E[A|Z=1,L]-E[A|Z=-1,L]$ almost surely.
 
\item Assumption~A of \cite{Han2020}: The following two conditions hold given $L$,\\ (a) either $E[Y_1|L, U] \geq  
E[Y_{-1}|L, U]$ or $E[Y_1|L, U] \leq E[Y_{-1}|L, U]$ almost surely; and \\
(b) either $E[A_1|L, U] \geq  
E[A_{-1}|L, U]$ or $E[A_1|L, U] \leq E[A_{-1}|L, U]$ almost surely.
\end{itemize}

\bibliographystyle{asa}
\bibliography{causal,iv}

\begin{thebibliography}{23}
\newcommand{\enquote}[1]{``#1''}
\expandafter\ifx\csname natexlab\endcsname\relax\def\natexlab#1{#1}\fi

\bibitem[{Angrist et~al.(1996)Angrist, Imbens, and Rubin}]{angrist1996}
Angrist, J.~D., Imbens, G.~W., and Rubin, D.~B. (1996), \enquote{Identification
  of Causal Effects Using Instrumental Variables,} \textit{Journal of the
  American Statistical Association}, 91, 444--455.

\bibitem[{Balke and Pearl(1997)}]{Balke1997}
Balke, A. and Pearl, J. (1997), \enquote{Bounds on Treatment Effects from
  Studies with Imperfect Compliance,} \textit{Journal of the American
  Statistical Association}, 92, 1171--1176.

\bibitem[{Chakraborty and Moodie(2013)}]{chakraborty2013statistical}
Chakraborty, B. and Moodie, E. (2013), \textit{Statistical methods for dynamic
  treatment regimes}, Springer.

\bibitem[{Cui and Tchetgen~Tchetgen(2020)}]{cui2020}
Cui, Y. and Tchetgen~Tchetgen, E. (2020), \enquote{A Semiparametric
  Instrumental Variable Approach to Optimal Treatment Regimes Under
  Endogeneity,} \textit{Journal of the American Statistical Association}, 0,
  1--12.

\bibitem[{Greenland(2000)}]{10.1093/ije/29.4.722}
Greenland, S. (2000), \enquote{{An introduction to instrumental variables for
  epidemiologists},} \textit{International Journal of Epidemiology}, 29,
  722--729.

\bibitem[{Han(2019)}]{han2019optimal}
Han, S. (2019), \enquote{Optimal Dynamic Treatment Regimes and Partial Welfare
  Ordering,} \textit{arXiv preprint arXiv:1912.10014}.

\bibitem[{Han(2020)}]{Han2020}
--- (2020), \enquote{Comment: Individualized Treatment Rules Under
  Endogeneity,} \textit{Journal of the American Statistical Association}.

\bibitem[{Hernan and Robins(2006)}]{hernan2006epi}
Hernan, M. and Robins, J. (2006), \enquote{Instruments for Causal Inference: An
  Epidemiologist's Dream?} \textit{Epidemiology (Cambridge, Mass.)}, 17,
  360--72.

\bibitem[{Kallus et~al.(2019)Kallus, Mao, and Zhou}]{pmlr-v89-kallus19a}
Kallus, N., Mao, X., and Zhou, A. (2019), \enquote{Interval Estimation of
  Individual-Level Causal Effects Under Unobserved Confounding,} in
  \textit{Proceedings of Machine Learning Research}, eds. Chaudhuri, K. and
  Sugiyama, M., PMLR, vol.~89 of \textit{Proceedings of Machine Learning
  Research}, pp. 2281--2290.

\bibitem[{Kallus and Zhou(2018)}]{kallus2018confounding}
Kallus, N. and Zhou, A. (2018), \enquote{Confounding-robust policy
  improvement,} in \textit{Advances in neural information processing systems},
  pp. 9269--9279.

\bibitem[{Kosorok and Laber(2019)}]{kosorok2019review}
Kosorok, M.~R. and Laber, E.~B. (2019), \enquote{Precision Medicine,}
  \textit{Annual Review of Statistics and Its Application}, 6, 263--286.

\bibitem[{Murphy(2003)}]{murphy2003optimal}
Murphy, S.~A. (2003), \enquote{Optimal dynamic treatment regimes,}
  \textit{Journal of the Royal Statistical Society: Series B (Statistical
  Methodology)}, 65, 331--355.

\bibitem[{Qian and Murphy(2011)}]{qian2011performance}
Qian, M. and Murphy, S.~A. (2011), \enquote{Performance guarantees for
  individualized treatment rules,} \textit{Annals of statistics}, 39, 1180.

\bibitem[{Qiu et~al.(2020)Qiu, Carone, Sadikova, Petukhova, Kessler, and
  Luedtke}]{qiu2020}
Qiu, H., Carone, M., Sadikova, E., Petukhova, M., Kessler, R.~C., and Luedtke,
  A. (2020), \enquote{Optimal Individualized Decision Rules Using Instrumental
  Variable Methods,} \textit{Journal of the American Statistical Association},
  0, 1--18.

\bibitem[{Robins et~al.(2000)Robins, Hern{\'a}n, and
  Brumback}]{robins2000marginalSM}
Robins, J.~M., Hern{\'a}n, M.~A., and Brumback, B.~A. (2000), \enquote{Marginal
  structural models and causal inference in epidemiology.}
  \textit{Epidemiology}, 11 5, 550--60.

\bibitem[{Rubin and van~der Laan(2012)}]{rubin2012statistical}
Rubin, D.~B. and van~der Laan, M.~J. (2012), \enquote{Statistical issues and
  limitations in personalized medicine research with clinical trials,}
  \textit{The international journal of biostatistics}, 8, 18.

\bibitem[{Tsiatis et~al.(2019)Tsiatis, Davidian, Holloway, and
  Laber}]{tsiatis2019dynamic}
Tsiatis, A.~A., Davidian, M., Holloway, S.~T., and Laber, E.~B. (2019),
  \textit{Dynamic Treatment Regimes: Statistical Methods for Precision
  Medicine}, CRC Press.

\bibitem[{Wald(1940)}]{wald1940fitting}
Wald, A. (1940), \enquote{The fitting of straight lines if both variables are
  subject to error,} \textit{The annals of mathematical statistics}, 11,
  284--300.

\bibitem[{Wang and Tchetgen~Tchetgen(2018)}]{wang2018bounded}
Wang, L. and Tchetgen~Tchetgen, E. (2018), \enquote{Bounded, efficient and
  multiply robust estimation of average treatment effects using instrumental
  variables,} \textit{Journal of the Royal Statistical Society: Series B
  (Statistical Methodology)}, 80, 531--550.

\bibitem[{Yadlowsky et~al.(2018)Yadlowsky, Namkoong, Basu, Duchi, and
  Tian}]{yadlowsky2018bounds}
Yadlowsky, S., Namkoong, H., Basu, S., Duchi, J., and Tian, L. (2018),
  \enquote{Bounds on the conditional and average treatment effect with
  unobserved confounding factors,} \textit{arXiv preprint arXiv:1808.09521}.

\bibitem[{Zhang et~al.(2012{\natexlab{a}})Zhang, Tsiatis, Davidian, Zhang, and
  Laber}]{zhang2012estimating}
Zhang, B., Tsiatis, A.~A., Davidian, M., Zhang, M., and Laber, E.
  (2012{\natexlab{a}}), \enquote{Estimating optimal treatment regimes from a
  classification perspective,} \textit{Stat}, 1, 103--114.

\bibitem[{Zhang et~al.(2012{\natexlab{b}})Zhang, Tsiatis, Laber, and
  Davidian}]{zhang2012robust}
Zhang, B., Tsiatis, A.~A., Laber, E.~B., and Davidian, M. (2012{\natexlab{b}}),
  \enquote{A robust method for estimating optimal treatment regimes,}
  \textit{Biometrics}, 68, 1010--1018.

\bibitem[{Zhao et~al.(2012)Zhao, Zeng, Rush, and Kosorok}]{zhao2012estimating}
Zhao, Y., Zeng, D., Rush, A.~J., and Kosorok, M.~R. (2012), \enquote{Estimating
  individualized treatment rules using outcome weighted learning,}
  \textit{Journal of the American Statistical Association}, 107, 1106--1118.

\end{thebibliography}

\end{document}